\documentstyle{amsppt}
\magnification=\magstep 1
\hsize=30truecc
\baselineskip=16truept
\TagsOnRight
\NoBlackBoxes
\leftheadtext{V\. Pestov}
\rightheadtext{Two questions of Mazur}
\newcount\refAue
\newcount\refBS
\newcount\refMau
\newcount\refno
\refno=0
\advance\refno by 1\refAue=\refno
\advance\refno by 1\refBS=\refno
\advance\refno by 1\refMau=\refno
\def\norm #1{{\left\Vert\,#1\,\right\Vert}}
\def\R {{\Bbb R}}
\def\C{{\Bbb C}}
\def\K{{\Bbb K}}
\def\e{{\varepsilon}}
\def\N{{\Bbb N}}
\rightline{\eightit Research Report RP-97-223 (May 1997),}
\rightline{\eightit School of Mathematical and Computing Sciences,}
\rightline{\eightit Victoria University of Wellington}
\vskip 1cm
\topmatter 
\title
Two 1935 questions of Mazur \\ about polynomials
in Banach spaces: \\ a counter-example
\endtitle 
\author 
Vladimir Pestov
\endauthor
\affil 
Victoria University of Wellington, New Zealand 
\endaffil
\address
School of Mathematical and Computing Sciences,
Victoria University of Wellington, P.O. Box 600, Wellington,
New Zealand.
\endaddress
\email vladimir.pestov$\@$mcs.vuw.ac.nz
\endemail

\abstract{We construct a continuous scalar-valued $2$-polynomial, $W$,
on the separable Hilbert space $l_2$
and an unbounded set $R\subset l_2$ such that (i) $W$ is bounded
on an $\varepsilon$-neighbourhood of $R$; (ii) $W$ is unbounded on
$\frac 1 2 R$; (iii) consequently, $W$ does not factor through any bounded
$1$-polynomial on $l_2$ sending $R$ to a bounded set.
This answers in the negative two 1935 questions 
asked by Mazur (problems 55 and 75 in the Scottish Book). 
The construction is valid both over $\R$ and $\C$.
(In finite dimensions the questions
were answered in the positive by Auerbach soon after being asked.)}

\endabstract
\endsubjclass
\endkeywords

\endtopmatter 
\document 

\subheading{1}
The following two problems were entered up in the Scottish Book by Mazur
in 1935. (We quote Ulam's translation from the 1981 edition
\cite{\the\refMau}.)

\proclaim{55 (Mazur)} There is given, in an $n$-dimensional space
$E$ or, more generally, in a space of type 
\rom($B$\rom), a polynomial $W(x)$
bounded in an $\e$-neighbourhood of a certain nonbounded set $R\subset E$
\rom(an $\e$-neighbourhood of a set $R$ is the set of all points which are
distant by less than $\e$ from $R$\rom). Does there exist a polynomial
$V(x)$ and a polynomial of first degree $\phi(x)$ such that
\smallskip
\item{\rm (1)} $W(x) = V(\phi(x))$;
\item{\rm (2)} The set $\phi(R)$, that it to say the image of the set
$R$ under the mapping $\phi(x)$, is bounded?
\endproclaim

\proclaim{75 (Mazur)} In the Euclidean $n$-dimensional space $E$ or,
more generally, in a space of type \rom($B$\rom) 
there  is given a polynomial
$W(x)$.
$\alpha$ is a number $\neq 0$. If a polynomial $W(x)$ is bounded in
an $\e$-neighbourhood of a certain set $R\subset E$ is it then bounded
in a $\delta$-neighbourhood of the set $\alpha R$
\rom(which is the set composed of elements $\alpha x$ for
$x\in R$\rom)? \rom(See problem 55.\rom)
\endproclaim

Clearly, an affirmative answer to problem 55 implies
an affirmative answer to problem 75. Auerbach wrote addenda to both
problems ({\it ibid.}) informing that he solved in the positive
problem 55 and, {\it ipso facto,} problem 75 in the
case of a (both real and complex) finite-dimensional space $E$ 
\cite{\the\refAue}.  The general case where $E$ is
a Banach space seems to have been open ever since.

In this note we answer the two questions in the negative by presenting
a counter-example where $W$ is a $2$-polynomial from $l_\infty$ to
itself. It can be easily converted into a scalar-valued $2$-polynomial on a
separable Hilbert space, and thus the answer remains in general
negative even in
the  most favourable infinite-dimensional setting. The construction is
the same in either real or complex case.

\subheading{2}
Let $\K=\R$ or $\C$.
For each $n\in\N$, define homogeneous polynomials of degree
$1$ and $2$, respectively:
$$W_n^{(1)}(x)=-\frac {nx}{n+1},$$
$$W_n^{(2)}(x)= \frac{x^2}{n+1}.$$
Since for each $n\in\N$ and 
$x\in\K$ one has $\vert W_n^{(1)}(x)\vert\leq\vert
x\vert$ and $\vert W_n^{(2)}(x)\vert\leq\vert x\vert^2$, the formula
$$l_\infty\ni x\equiv (x_n)_{n\in\N} \mapsto
\left(W_n^{(i)}(x_n)\right)_{n\in\N}, ~~ i=1,2$$
correctly determines homogeneous polynomials 
$W^{(1)},W^{(2)}\colon l_\infty\to l_\infty$ of degree $1$ and $2$,
respectively. Both $W^{(1)}$ and $W^{(2)}$ are continuous, being
bounded on the unit ball $B_1(0)_{l_\infty}$.
(Cf. \cite{\the\refBS}, Th. 1 and a remark on p. 61.)
Define a continuous $2$-polynomial mapping $W\colon l_\infty\to l_\infty$
via
$$W = W^{(1)}+W^{(2)},$$
that is,
$$l_\infty\ni x\equiv (x_n)_{n\in\N} \overset W\to\mapsto \left(
W_n(x_n)\right)_{n\in\N}\in l_\infty,$$ 
where for all $n\in\N$ and $x\in\K$
$$W_n(x)=W_n^{(1)}(x)+W_n^{(2)}(x)\equiv \frac{x(x-n)}{n+1}.$$
The following is obvious.
\proclaim{Claim 1} If either $\vert x\vert\leq 1$ or $\vert x-n\vert\leq 1$,
then $\vert W_n(x)\vert\leq 1$.
\endproclaim

Let $e_n$ is the `$n$-th basic vector' in $l_\infty$,
$(e_n)_k=\delta_{n,k}$. 
Define an unbounded subset 
$$R=\{ne_n\colon n\in\N\}\subset l_\infty.$$
\proclaim{Claim 2}
The polynomial $W$ is bounded on the $1$-neighbourhood of
$R$.
\endproclaim
$\triangleleft$ Let $x\in B_1(R)$, that is, for some $k\in\N$,
$x\in B_1(ke_k)$. It means 
$\vert x_k-k\vert\leq 1$ and $\vert x_n\vert\leq 1$
for all $n\neq k$. By Claim 1, $\vert W_n(x_n)\vert\leq 1$
for all $n$, that is,
$\norm{W(x)} \leq 1$. $\triangleright$
\proclaim{Claim 3} The polynomial $W$ is unbounded on $(1/2)R$.
\endproclaim
$\triangleleft$ 
For each $k\in\N$, one has 
$(k/2)e_k\in (1/2)R$ and
$$\norm{W\left(\frac k 2 e_k\right)}=\left\vert W_k\left(\frac k 2
\right)\right\vert=\frac{k^2}{2(k+2)}\to\infty~~
\text{as}~~k\to\infty.\phantom{xxx}\triangleright$$

\subheading{3} It is easy to further modify the example so as to
obtain a $\K$-valued $2$-polynomial on $l_2$ in possession of the same
properties. Applying the Banach--Steinhaus theorem to the
unbounded set $W((1/2)R)\subset l_\infty$, choose a bounded
linear functional
$\varphi\colon l_\infty\to\K$ with
$\varphi[W((1/2)R)]$ unbounded. Denote by
$\widetilde W$ the  
composition of the three mappings:
$$l_2\overset i\to\hookrightarrow l_\infty\overset{W}\to\to l_\infty
\overset{\varphi}\to\to\K,$$
where $i\colon l_2\hookrightarrow l_\infty$ stands for the canonical
contractive injection. The mapping
$\widetilde W\colon l_2\to\K$ is 
continuous $2$-polynomial.
Set $\widetilde R=i^{-1}(R)$. Since 
$i(B_1(R)_{l_2})\subset B_1(R)_{l_\infty}$, the polynomial
$\widetilde W$ is  bounded on the
$1$-neighbour\-hood of $\widetilde R$, and since  
$i(\widetilde R)=R$, it is unbounded on $(1/2)\widetilde R$.

\enddocument

\Refs 
\widestnumber\key{OdlH22}
\vskip0.3truecm

\ref\key\the\refAue
\by H\. Auerbach
\paper Un th\'eor\`eme sur les polynomes \`a $n$ variables
\jour Studia Math\.
\vol 5
\yr 1934
\pages 171--173
\endref

\ref\key\the\refBS
\by J\. Bochnak and J\. Siciak
\paper Polynomials and multilinear mappings in topological vector
spaces
\jour Studia Math\.
\vol 39
\yr 1971
\pages 59--76
\endref

\ref\key\the\refMau
\by R\.D\. Mauldin (ed.)
\book The Scottish Book
\publ Birkh\"auser
\publaddr Boston--Basel--Stuttgart
\yr 1981
\endref


\endRefs
\bye